\documentclass[10pt]{article}
\usepackage[letterpaper,top=47.2mm,bottom=47.2mm,left=45.45mm,right=45.45mm,includefoot]{geometry}
\usepackage{verbatim}
\usepackage{amsfonts} 
\usepackage{graphicx}
\usepackage{amssymb}
\usepackage{epsfig}

\newtheorem{theorem}{Theorem}

\newtheorem{lemma}{Lemma}
\input epsf
\begin{document}
\title{The Quandle and Group for General Pairs of Spaces\footnote{MSC 2000: 17D99, 57R40.}}
\author{Blake Winter\\
\emph{University at Buffalo}\\
\emph{Email: bkwinter@buffalo.edu}}
\maketitle
\begin{abstract}
Joyce has shown that the fundamental quandle of a classical
knot can be derived from consideration of the fundamental
group and the peripheral structure of the knot, and also
that the group and much of the peripheral structure can be recovered from the
quandle. We generalize these results to arbitrary dimensions,
and also to virtual and welded knots and arcs.\end{abstract}
\section{Quandles}
\begin{figure}[htbp]
\begin{centering}
\includegraphics[scale=1]{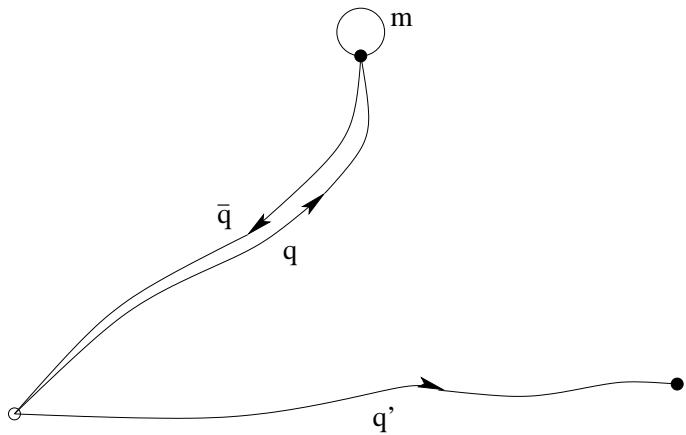}
\caption{The quandle operation for the fundamental quandle of a pair of spaces. The
open circle indicates the basepoint of $M$, while the black circles represent
points on $\partial \overline{N}$.\label{qop}}
\end{centering}
\end{figure}

A \emph{quandle} (sometimes called a \emph{distributive groupoid})\cite{DJ, Mat}
is a non-associative algebraic structure which is particularly connected
to the study of classical knots. In particular a quandle may be defined as a set $Q$
together with an left-invertible binary operation which we will write by exponentiation
$a^{b}$ obeying the following relations:
\begin{equation}
  a^{a}=a\label{q1}
\end{equation}
\begin{equation}
	(a^{b})^{c}=(a^{b})^{b^{c}}\label{q2}
\end{equation}
Note that left invertibility is equivalent to also requiring that the equation $x^{a}=b$
with variable $x$ should have exactly one solution. These relations are analogues of
the Reidemeister moves for classical knots. In particular, Eq. \ref{q1} reflects
the first Reidemeister move, left invertibility reflects the second, and
Eq. \ref{q2} reflects the third. For this reason, they are particularly suited to the study of knots. Joyce\cite{DJ}
showed in particular that there is a quandle naturally associated to any knot
in $S^{3}$ which is a nearly complete (up to mirror-reversal) invariant of the knot. The quandle
associated to a classical knot is a special case of Joyce's quandle associated to any codimension two knotting. We briefly
review this construction here, referring the reader to Joyce's exposition for
further details.\\
Let $M$ be a path-connected manifold, and let $K$ be a codimension two knotting
(neither $M$ nor $K$ need to be spheres).
Let $N$ be a tubular neighborhood of $K$. Then $M-N$ is a manifold with boundary.
One component of $\partial M$ will be $\partial \overline{N}$. Let $P$ be
the subgroup of $\pi_{1}(M-N)$ which is the image of the homomorphism induced
by the injection $\partial \overline{N}\rightarrow M-N$. We refer to this
as the \emph{peripheral} subgroup. Note that $P$ is defined only up to conjugation;
any of its conjugates are also peripheral.\\
Now suppose that in $P$ there is a preferred element $m$,
the \emph{meridian}. In particular for any codimension two knot, the meridian
is canonically defined up to conjugation. 
We may now define the quandle $Q$ associated to the pair $(M,N)$ and the chosen
meridian. The elements of $Q$ are homotopy classes of paths that start at
the basepoint of $M$ (that is, the basepoint chosen for the fundamental group)
and end in $\partial \overline{N}$, where homotopies are required to preserve
these two conditions. The quandle operation on two elements of $Q(M-N)$ is defined
geometrically as illustrated in Fig. \ref{qop}. In particular, $[q]^{[q']}$
is the homotopy class of paths containing the path $q'\overline{q}mq$, where
the overbar indicates that the path is followed in reverse. It is
a striaghtforward exercise to show that this is well-defined and meets
the definition of a quandle operation.
Note that there is an action $\pi_{1}(M-N)$ action on $Q(M-N)$: $[q]g$ is the equivalence
class including the homotopy class of the path $qg$. We follow the notation
of Joyce for this action.
\section{Quandles From Groups}
Given a group $G$, a subgroup $P$, and an element $m$ in the center $Z(P)$, Joyce
defines a quandle as follows. The underlying set is the set of right cosets of $P$, $P\backslash G$.
We define the quandle operation such that $Pg^{Ph}=P(gh^{-1}mh)$. This is well-defined,
because if $h'\in Ph$, then $h'=ah$ for some $a\in P$. Then
\begin{equation}
  h'^{-1}mh'=h^{-1}a^{-1}mah=h^{-1}mh,
\end{equation}
where the last equality follows from $m\in Z(P)$.
We will denote a quandle constructed in this way
by $(P\backslash G,m)$, or simply $P\backslash G$ when there is a canonical choice for $m$ (as there is
when $m$ is a meridian and $P$ a peripheral subgroup).
\section{Group Actions on Quandles}
\begin{figure}[htbp]
\begin{centering}
\includegraphics[scale=1]{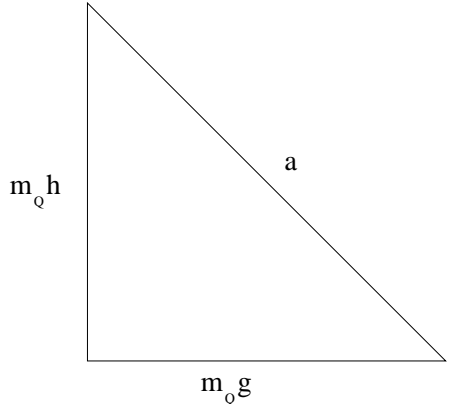}
\caption{There is a homotopy of $m_{Q}g$ to $m_{Q}h$ with the endpoint of the
path tracing out a path $a$ on $\partial N$.\label{tri}}
\end{centering}
\end{figure}

The work in this section is largely a generalization of the work of Joyce\cite{DJ} for
the case of classical knots.\\
Consider a codimension two pair $(M,K)$, and the tubular neighborhood $N$
of $K$. Then we can define the quandle $Q=Q(M,K)$
geometrically as in the first section. On the other hand we can define the quandle
$G$ to be the quandle $(P\backslash \pi_{1}(M-K),m)$ for P some arbitrarily chosen peripheral
subgroup of $\pi _{1}(M-K)$, and $m$ the meridian in $P$.\\
Now $\pi _{1}$ acts on $Q$ by setting $[q]g$ to be the equivalence
class including the homotopy class of the path $qg$. $m$ is an element of
the fundamental group, but we may also define an analogous element $m_{Q}$ in
the fundamental quandle to be the quandle element such that $\overline{q}mq$
is the meridian of $P$.
\begin{lemma}
Every element $q\in Q$ may be obtained as $m_{Q}g$
for some element $g\in \pi_{1}(M-N)$.\label{transq}
\end{lemma}
\textbf{Proof:} We choose representatives of $q,m_{Q}$ whose endpoints on $\partial N$
agree. Then the path $m_{Q}\overline{m_{Q}}q$ represents the same quandle class as $q$.
However, the path $\overline{m_{Q}}q$ starts and ends at the basepoint, and hence represents
some element of the fundamental group.
Thus (abusing the notational difference between paths and their
equivalence classes) $q=m_{Q}\overline{m_{Q}}q=m_{Q}g$.$\square$\\
This shows in addition that the group action on $Q$ is transitive.
\begin{lemma} $m_{Q}g=m_{q}h$ iff $h=ag$ for some $a\in P$.\end{lemma}
\textbf{Proof:} If $h=ag$ then it is straightforward to see that $m_{Q}g=m_{q}h$.
On the other hand if $m_{Q}g=m_{q}h$, then there is a map of the disk into $M$
as shown in Fig. \ref{tri}.
But this disk shows that the path which follows $\overline{h}\overline{m_{Q}}am_{Q}g$
is contractible. However, $\overline{m_{Q}}am_{Q}$ is in the peripheral group,
and we have our result.$\square$\\
A similar argument shows that the stabilizer of $m_{Q}g$ is $g^{-1}Pg$. For elements
in that conjugate peripheral group are easily seen to stabilize $m_{Q}g$. In
particular if $p\in P$, then $m_{Q}g(g^{-1}pg)=m_{Q}pg=m_{Q}g$. On the other hand
if some some element of the fundamental group stabilizes $m_{Q}g$ we may build
another homotopy disk as in the previous lemma. As the stabilizer does not play
a role in our following argument we leave the details of the construction as an exercise.\\
Let us now turn our attention to the quandle $G$ which is built from the fundamental
group and the peripheral information. Our goal is to establish the following result:
\begin{theorem}\label{main}
$Q\cong G$ as quandles; hence the quandle can be recovered from the fundamental
group and peripheral information.
\end{theorem}
We begin by noting that the underlying set of $G$ is just the set of right cosets
of $\pi_{1}(M-K)$. Therefore $\pi_{1}(M-K)$ naturally acts on $G$, and this action
is obviously transitive. Indeed, every element of $G$ is of the form
$Pg$ for some $g\in \pi_{1}(K)$. Furthermore, $Pg=Ph$ iff $h=ag$ for some
$a\in P$.\\
We now define a quandle morphism $Q\rightarrow G$ sending $m_{Q}g\mapsto Pg$.
This is well defined and invertible by our previous lemmas, so we need only
show that it is a quandle map as well. Now 
\begin{equation}
  (m_{Q}g)^m_{Q}h=m_{Q}g\overline{h}\overline{m_{Q}}mm_{Q}h=m_{Q}g\overline{h}mh.
\end{equation}
This is mapped to 
\begin{equation}
  Pg\overline{h}mh,
\end{equation}
which is precisely $Pg^{Ph}$. That the inverse map is also
a morphism of quandles follows similarly. Therefore this is an invertible quandle morphism,
and hence an isomorphism of quandles, establishing Thm. \ref{main}.$\square$
\section{Remarks}
Joyce has also proved a theorem for classical knots
stating that the triple $(\pi_{1}(M-K), P, m)$ can be reconstructed from $Q(M-K)$. This
is because $\pi_{1}(M-K)$ is isomorphic to $Adconj(Q)$, the free group on $Q(M-K)$ quotiented
by taking all the relations on $Q(M-K)$ and treating the quandle operation as conjugation in the group.
We write an element $q\in Q$ as $\hat{q}$ when we wish to designate them as elements of $Adconj(Q)$.
Note that $Adconj(Q)$ is canonically isomorphic to $\pi_{1}(M-K)$ (up to a choice of basepoint),
by sending $\hat{q}$ to the homotopy class of a Wirtinger generator, $[\overline{q}m_{Q}q]$.
Then in $Adconj(Q)$ one may take $m$ to be any of the elements of the generators defined. $Adconj(Q)$ acts
on $Q$ by $q\hat{q_{0}}...\hat{q_{n}}=(...(q^{q_{0}})^{q_{1}}...)^{q_{n}}$. If an element
of $Adconj(Q)$ has two presentations, $\hat{q_{0}}...\hat{q_{n}}$ and $\hat{p_{0}}...\hat{p_{m}}$,
this implies that those products of Wirtinger generators are homotopic rel basepoints, and
this homotopy passes to a homotopy of $(...(q^{q_{0}})^{q_{1}}...)^{q_{n}}$ to $(...(q^{p_{0}})^{p_{1}}...)^{p_{m}}$.
It is then straightforward to check that this action of $Adconj(Q)$ is just the geometrically
defined action of $\pi_{1}(M-K)$ on $Q$.
Then $P$
will be the stabilizer of $m\in Q$ under this action, and so the triple can be recovered
from the quandle.\\
However, this only permits us to recover
the triple when all the relations on $\pi_{1}(M)$ can be written down in terms of
conjugation. If $\pi_{1}(M)$ has no such presentation, then this will not be possible.
Checking this condition is complicated, however, by the fact that $Q$ will generally
be infinite, and so the presentation we obtain will not be a finite presentation. Nonetheless,
for higher-dimensional knots in spheres, the above argument shows:
\begin{theorem}
For codimension two knots in spheres of any dimension, the quandle and the triple $(\pi_{1}(M-K), P, m)$
can both be calculated from the other. More generally this holds for codimension two knots
in simply connected spaces.\label{hid}\end{theorem}
Thm. \ref{main} does imply that any invariant based upon quandles for any knot should
be interpretable as an invariant based upon classical knot invariants. Thus,
for example, the quandle cocycle invariants\cite{CKS} are determined by the triple $(\pi_{1}(M-K), P, m)$,
even for higher-dimensional knots. This becomes even more powerful when it is recalled
that the peripheral group of a sphere is cyclic, and hence is determined by $m$. On the other
hand, since in general $P\cong(\pi_{1}(K)/L)\times \mathbb{Z}$ for some normal subgroup $L$,
this indicates that for higher-dimensional knots which as manifolds have large
fundamental groups, the quandle has the potential to capture more information
than for those with smaller fundamental groups.\\
Eisermann\cite{ME} has in fact shown that 
the quandle cocycle invariants of classical knots is a specialization of certain colourings of their fundamental
groups. His construction makes use of the full peripheral structure of the classical knot
(that is, the longitude and the meridian),
and hence does not generalize immediately to higher dimensions. However, in light
of our result here, we pose the question of whether a similar construction might
not be possible in higher dimensions.\\
As a last application of our result, we give the following theorem.
Recall that \emph{virtual} knots are a combinatorial generalization of knot diagrams
introduced by Kauffman\cite{Kauf}, and welded knots are a quotient of virtual
knots given by adding an additional permitted move, first explored for braids
in \cite{FRR}. When the diagrams are allowed to have two endpoints instead of
being closed we obtain virtual (welded) \emph{arcs}\cite{SS,CKS}.
\begin{theorem}
The fundamental quandle of a virtual or welded knot (or arc) $K$ determines, and is determined by
the triple $(\pi_{1}(K), P, m)$.\end{theorem}
\textbf{Proof:} Satoh\cite{SS} has defined a map $Tube$ on welded knots and arcs (hence on virtual knots and arcs)
which maps each knot or arc to a surface knot, and shown that $Tube$ preserves the quandle and fundamental group.
In addition, a direct calculation shows
that $P$ and $m$ are preserved by $Tube$\cite{BW}.
But the latter determine the quandle for the surface knot, and hence also for the virtual or welded knot or arc.\\
To see that the quandle determines the triple $(\pi_{1}(K), P, m)$, it suffices to check this
for surface knots, as was done in Thm. \ref{hid}.$\square$

\end{document}